%Paper: q-alg/9502016
%From: Murray Gerstenhaber <murray@math.upenn.edu>
%Date: Thu, 23 Feb 95 13:38:38 -0500

\input amstex
\documentstyle{amsppt}
\leftheadtext{Murray Gerstenhaber and Mary E. Schaps}
\rightheadtext{Hecke algebras and the Donald--Flanigan conjecture
for $S_n$}

\magnification=1200
\parindent 20pt
\NoBlackBoxes
\topmatter
\title Hecke algebras, $U_qsl_n$, and the Donald--Flanigan
conjecture for $S_n$
\endtitle
\author Murray Gerstenhaber and Mary E. Schaps
      \endauthor
\address Dept. of Mathematics, University of Pennsylvania,
Philadelphia, PA 19104-6395
      \endaddress
\email mgersten\@mail.sas.upenn.edu or murray\@math.upenn.edu
 \endemail
\address Dept. of Mathematics and Computer Science, Bar Ilan
University, Ramat-Gan 52900, Israel
      \endaddress
\email mschaps\@bimacs.cs.biu.ac.il
      \endemail
\date October 4, 1994 \enddate
\thanks Research of the first author was partially supported by a
grant from the NSA. \endthanks

\abstract
To each partition $\frak p$ of $n$ we associate in a canonical
way a simple $S_n$ module with an orthogonal basis indexed by
Young diagrams in a way which carries over immediately to the
quantized case.  With this we show that the Hecke algebra of
$S_n$ is a global solution to the Donald--Flanigan problem for
$S_n.$  The procedure gives ``canonical'' primitive idempotents
different from the classical ones of Frobenius--Young and makes
some number--theoretic statements.
\endabstract
\endtopmatter

\document
\baselineskip 12pt

\subhead 1. Introduction \endsubhead
The main result of this paper is a new demonstration of the
relation between (non-increasing) partitions $\frak p$ of $n$ and
representations of the symmetric group $S_n.$  Our method
carries over to the ``quantized'' case, where the group algebra
of $S_n$ is replaced by its Hecke algebra $\Cal H_n,$  and yields
the following:  To each $\frak p$ there is canonically associated
a simple module with a non-degenerate, symmetric bilinear form,
for brevity called an inner product, and a canonical orthogonal
basis indexed by the various Young diagrams associated to $\frak
p;$  the rank of this module is the number of such diagrams. That
these are all the simple modules and that the Hecke algebra
operates as the full ring of linear endomorphisms of each will be
evident.  In the unquantized case our procedure produces
primitive idempotents which, although also indexed by Young
diagrams, are different from the classical ones computed using
Young symmetrizers. An incidental result  is certain
number--theoretic statements which we do not know how to
interpret.

Of greatest importance to us for application to the Donald--
Flanigan conjecture is the ring of coefficients over which all
this takes place.  The group ring of a finite group is separable
over any ring in which its order is invertible. For $S_n$ it is
therefore sufficient to take as coefficients $\Bbb Z_{1,n} :=
\Bbb Z[1/n!].$ Then $\Bbb Z_{1,n}S_n$ is not only separable but a
direct sum of total matric algebras over $\Bbb Z_{1,n}.$  Set
$\Bbb
Z_q := \Bbb Z[q,q^{-1}]$ where $q$ is a variable. The {\bf Hecke
algebra} $\Cal H_n$ of $S_n$ is a free module of rank $n!$ over
$\Bbb Z_q$ with basis elements $T_w$ indexed by the elements
$w\in  S_n$ and multiplication given as follows: The {\it length}
$\ell(w)$ is the number of factors in a shortest expression of
$w$ as a product of generators $s_i := (i,i+1), i = 1,\dots,n-1,$
of $S_n.$  Multiplication is determined by setting (i) $T_sT_w =
T_{sw}$ if $s$ is one of these generators and  $w\in S_n$ is an
element with
$\ell(sw) > \ell(w)$, and (ii) $T_s^2 = (q-q^{-1})T_s+1$.  This
implies that $T_sT_w = (q-q^{-1})T_w + T_{sw}$ when $\ell(sw) <
\ell(w)$.  (Often instead of (ii) one takes $T_s^2 = (q-1)T_s +
q$, cf. \cite{Hu}; that definition can be transformed into ours
by substituting $q^2$ for $q$ and dividing the generators by the
new $q.$  The present form is more useful when dealing with
quantum groups.)  When necessary to indicate the dependence on
the
parameter $q$ we may write $\Cal H_n(q).$  Writing $1+t$ for $q$
one sees that $\Cal H_n(1+t)$ is a deformation of $\Bbb ZS_n.$

Set $i_q:= (1-q^i)/(1-q)$ and similarly $i_{q^2}:=
(1-q^{2i})/(1-q^2).$ These are the ``$q$--numbers''. For $i$ an
integer they are polynomials in $q$ (or $q^{-1}$ for $ i < 0$.)
We call $i$ the ``argument '' and $q$ the ``parameter'' of $i_q.$
Set $n_{q^2}! := n_{q^2}(n-1)_{q^2}\dots 2_{q^2}$ and $\Bbb
Z_{q,n} = \Bbb Z[q, q^{-1}, 1/n_{q^2}!].$  We will prove the
following, from which one recovers the corresponding assertion
for $\Bbb Z[1/n!]S_n$ by letting $q \to 1.$

\proclaim{Theorem 1.1} Over $\Bbb Z_{q,n}$ the Hecke algebra
$\Cal H_n$ becomes a direct sum of total matric algebras.
\endproclaim

The original Donald--Flanigan conjecture asserts that if $G$ is a
group of finite order divisible by a prime $p$ then $\Bbb F_pG$
can be deformed to an $\Bbb F_p[[t]]$--algebra which becomes
separable when coefficients are extended to the Laurent series
field $\Bbb F_p((t)).$  Such an algebra will be called a {\it
solution} to the Donald--Flanigan problem for $G$ {\it at the
prime} $p.$  Extending coefficients to the algebraic closure of
$\Bbb F_p((t))$ gives an algebra which is a direct sum of total
matric algebras called its ``blocks''.  Different solutions to
the Donald--Flanigan problem at the same prime may have different
block sizes.  For example, if $G = \Bbb Z/2 \times \Bbb Z/2$ then
$\Bbb F_2G \cong \Bbb F_2[x]/(x^2) \otimes_{\Bbb F_2}\Bbb
F_2[y]/(y^2),$ which deforms to $\Bbb F_2[[t]][x]/(x^2+tx)
\otimes_{\Bbb F_2[[t]]}\Bbb F_2[y]/(y^2+ty).$  This commutative
algebra becomes separable when coefficients are extended to $\Bbb
F_2((t)),$ and its blocks are all one-dimensional.  These are the
same as for the complex group algebra $\Bbb CG.$  On the other
hand, there is also a non-commutative solution.  Let $\Bbb
F_2\langle x,y\rangle$  denote the
non--commutative polynomial ring in two variables over the field
of two elements.  Then $\Bbb F_2G \cong \Bbb F_2\langle
x,y\rangle/(x^2,y^2,xy-yx),$ which deforms to $\Bbb F_2\langle
x,y\rangle/(x^2, y^2,xy-yx-t).$ Over $\Bbb F_2((t))$ the latter
is just a $2 \times 2$ total matric algebra whose blocks do not
correspond to those of the previous solution.  (A solution in
which the block sizes correspond to those
over $\Bbb C$ is called {\it equimodular}.)  In view of this we
make the following
\proclaim{Definition} A global solution to the Donald--Flanigan
problem for a finite group $G$ is a deformation $A_t$ of the
integral group ring $\Bbb ZG,$ together with a multiplicatively
closed subset $S$ of the coefficient ring $\Bbb Z[[t]]$ such that
i) $S^{-1}A_t$ is separable over $S^{-1}\Bbb Z[[t]], $and ii) $S$
contains no rational prime dividing the order $\#G$ of $G.$
\endproclaim
If ii) is weakened to allow that $S$ contain certain rational
primes $p', p'', \dots$ dividing $\#G$ then we say that we have a
global solution away from $p',p'',\dots.$  When $S$ is
understood, we generally refer to the $\Bbb Z[[t]]$ algebra $A_t$
itself as the global solution.

The concept of a global solution is important even when $G$ is a
$p$--group. For notice that the composite map $\Bbb Z[[t]] \to
\Bbb F_p((t))$ consisting of reduction modulo $p$ and inversion
of $t$ (which may be done in either order) can be factored
through $S^{-1}\Bbb Z[[t]].$  Therefore, the separable algebras
$\Bbb F_p((t))\otimes_{S^{-1}\Bbb Z[[t]]}S^{-1}A_t$ for the
various primes $p$ dividing $\#G$ are all quotients of the same
separable algebra $S^{-1}\Bbb Z((t)) \otimes_{S^{-1}\Bbb
Z[[t]]}S^{-1}A_t$ and their blocks are therefore in natural
correspondence with the blocks of the latter. But the latter has
characteristic $0,$ so extending its coefficients to include
$\Bbb C$ we see that a global solution not only gives a ``local''
solution at each prime (or in the weaker case, away from primes
$p',p'',\dots$,) but naturally identifies the blocks of each of
the local solutions which it generates with those of the complex
group algebra $\Bbb CG.$ It is a scheme-theoretic solution
generalizing Maschke's theorem in the strongest possible way.
With this definition, Theorem 1.1 has the following immediate

\proclaim {Corollary 1.2} Setting $q = 1+t,$ the Hecke algebra
$\Cal H_n(1+t)$ together with the multiplicatively closed subset
of $\Bbb Z[[t]]$ generated by $1/n_{q^2}! = 1/n_{(1+t)^2}!$ is a
global solution to the Donald--Flanigan  problem for the
symmetric group $S_n.$
\endproclaim
Another immediate corollary to the first theorem is that $\Cal
H_n(q)$ splits (i.e., becomes a direct sum of total matric
algebras) over any field $k$ in which $q$ is not a $2i$--th root
of unity for any $i = 2,\dots n.$   One will find the proposition
in this form in Dipper and James \cite{DJ}. Previously Curtis,
Iwahori, and Kilmoyer \cite{CIK}, working over $\Bbb C,$ showed
that if $q$ is not a root of unity in the usual sense then the
Hecke algebra is isomorphic to $\Bbb CS_n.$ The hypothesis that
$k$ be a field is, however, too strong to permit one to conclude
that $\Cal H_n$ is a global solution to the Donald--Flanigan
conjecture. It does show that it is one at each individual prime
but gives no way of linking the primes to show that the matrix
blocks correspond. Wenzl \cite{W} proved a result closer to our
first theorem by showing explicitly, using ``quantized'' Young
symmetrizers, how to construct the idempotents of the localized
Hecke algebra away from the prime $2$.  That is, his method
requires extension of the coefficients to $\Bbb Z[q, q^{-1}, 1/2,
1/n_{q^2}!].$ (It fails at $2$ because Wenzl uses square roots,
which he says could be avoided.) This gives a global solution
away from $2,$ but in
group theory the prime $2$ is indispensable since by the
Feit--Thompson theorem a group of odd
order is solvable, \cite{FT}.

Rather than amend either \cite{W} or \cite{DJ}, we take an
approach which gives information about representations of
the symmetric group even over the rationals.  Our representations
are on the tensor powers $V^{\otimes n}$ of a free module $V,$
where the symmetric group is represented by permutations of the
tensor factors. In the ``quantized'' case, the generators of the
Hecke algebra are represented by quantum Yang--Baxter matrices.
These  generate the commutant of the standard quantized universal
enveloping algebra $U_qsl_n$ operating on the same space,
extending to the quantized case a central observation Schur's
thesis \cite{S}.  One immediately recovers, amongst others, the
basic result of Lusztig \cite{L} and Rosso \cite{R} that
$U_qsl_n$ has over $\Bbb C$ essentially the same representation
theory as $Usl_n$ provided that $q$ is not a root of unity. For
generic $q$ the cohomology theory of Hopf algebras \cite{GSk1,2;
Shn} gives the stronger result that any quantization, not merely
the standard one, in fact has this
property. For other relevant work, cf. also \cite{GGSk1, CFW, F},
as well as the thesis of P.N. Hoefsmit \cite{Hoe} which, although
often cited, unfortunately remains unpublished.

\subhead 2. Some elementary representation theory of $S_n$
\endsubhead
As an introduction to our approach, we review some of the
elementary representation theory of $S_n$. Let $V$ be a vector
space of dimension $d$ over $\Bbb Q$ with a basis $x_1,\dots,x_d$
except that when $d = 2$ we will write $x$ for $x_1$ and $y$ for
$x_2.$ The $n$--th tensor power of $V$ will be denoted simply
$V^n.$ On this $S_n$ operates by permutation of the tensor
factors, inducing an operation of its group algebra $\Bbb QS_n,$
and $\operatorname{End}V$ operates diagonally.  These operations
clearly commute and Schur \cite{S} shows that they are mutual
commutants inside $\operatorname{End} V^n,$ i.e., each consists
of all operators commuting with the other.  From a modern
viewpoint this is clear, since it is evident that any operator
commuting with $S_n$ must act diagonally, and since $\Bbb QS_n$
is separable its second commutant inside $\operatorname{End} V^n$
must be itself.

Tensor products of elements, like $x\otimes y,$ will be denoted
simply by concatenation, $xy.$ Note that if $d \ge n$ then the
representation of $\Bbb QS_n$ on $V^n$ is faithful, since the
$n!$ images of $x_1x_2\dots x_n$ are linearly independent.
Having chosen a basis for $V$ we can make it into an inner
product space by taking these basis elements to be orthonormal.
This induces an inner product $(a,b)$ on $V^n$ in which the
monomials of total degree $n$ in the basis elements of $V$ form
an orthonormal basis for $V^n.$  The adjoint of an $L \in
\operatorname{End} V^n$ will be denoted $L^t,$ so by definition
$(La,b) = (a,L^tb).$ Our inner product is symmetric, so the
matrix of $L^t$ relative to an orthonormal basis is just the
transpose of that of $L.$  The following is then obvious.

\proclaim{Lemma 2.1} If  $\sigma \in S_n$ then $\sigma^t =
\sigma^{-1}$; i.e., the elements of $S_n$ are orthogonal
operators. $\quad\square$
\endproclaim

Any real representation of a finite group is equivalent to one by
orthogonal matrices since we can introduce an invariant metric on
the representation space by first choosing an inner product
arbitrarily and then averaging with respect to the group
operations. The foregoing is, however, both canonical and extends
to the quantized case. It yields

\proclaim{Theorem 2.2} Central elements of $\Bbb QS_n$ are
self--adjoint. \endproclaim
\demo{Proof} The sum of all elements in any conjugacy class is
self--adjoint because in $S_n$ every element is conjugate to its
inverse, and the central elements are the linear combinations of
these sums. $\square$ \enddemo

Notice that the orthogonal projection of one $S_n$ submodule of
$V^n$ on another is a submodule of the latter which is a quotient
module of the former. We therefore have the following, which
extends to the quantized case.

\proclaim{Theorem 2.3} The isotypical components of $V^n$ are
mutually orthogonal. $\square$ \endproclaim

On the tensor algebra $TV$ of $V$, of which $V^n$ is a
homogeneous component, denote the derivation $\partial/\partial
x_i$ by $\partial_i.$  We then have a representation of the Lie
algebra $sl_n:$ letting $H_i, X_i, Y_i, i = 1, \dots,d-1$
be a Cartan basis, send $X_i$ to $x_i\partial_{i+1}, Y_i$ to
$x_{i+1}\partial_i$ and $H_i$ to their commutator
$x_i\partial_i - x_{i+1}\partial_{i+1}.$  These are the
``infinitesimal generators" of the special linear group
$Sl(V)$ acting diagonally, they commute with the action of
$S_n,$ and any linear operator commuting with all of them lies in
$\Bbb QS_n.$  (To see this, extend coefficients to $\Bbb R$ and
note that the special linear group is generated by the
exponentials of its infinitesimal generators.) The representation
of $sl_d$ gives rise to one of its universal enveloping
algebra $U:=Usl_d.$  For the proofs of the theorems
asserted in the Introduction it will never be necessary to use
the Casimir operator, but it is interesting to consider it. Since
any element of the center of $U$ is {\it a fortori} in its
centralizer, such an element must operate like a central element
of $\Bbb QS_n,$  so we have

\proclaim{Theorem 2.4} Central elements of $U$, and in particular
the quadratic and higher Casimir operators, act self-adjointly on
$V^n.$ Submodules of distinct eigenspaces of a central element of
$U$ are never isomorphic. $\square$ \endproclaim

For $sl_2$ the quadratic Casimir operator is $\Cal C =\frac 12
+\frac 12 H^2 + XY + YX.$  (The constant term $\frac 12$ is not
important for us but arises naturally and is essential in certain
parts of the theory of quantum groups when one must
exponentiate the Casimir element.)  It is easy to check in this
case that $\Cal C$ is the essentially unique quadratic central
element of $U.$ The eigenspaces of $\Cal C$ are $S_n$ submodules
of $V^n$ and distinct eigenspaces of a self--adjoint operator are
mutually orthogonal, so we have immediately a decomposition of
$V^n$ into an orthogonal direct sum of submodules.  This is
generally not  a full isotypical decomposition since the
quadratic Casimir of a simple Lie algebra $\frak g$ may have the
same eigenvalue on non--isomorphic simple $\frak g$ modules. The
center of the universal enveloping algebra of a simple Lie
algebra of rank $r$ is a polynomial ring in $r$ elements,
thequadratic and
higher'Casimir operators; together they give the full
isotypical decomposition. For $d=2$ the rank is $1$ and all
central elements are polynomials in the quadratic Casimir, so the
problem does not arise; it is trivial in this case $\Cal C$ has
distinct eigenvalues on non-isomorphic simple modules. (This is {\it a
fortiori} true in the quantized case; see the next section.)

Letting $d = 2,$ denote the basis elements of $V$ by $x, y$ and
those of the Lie algebra $sl_2$ by $H, X, Y,$ where $Xx =
0, Xy = x, Yx = y, Yy = 0, Hx = x,$ and $Hy = -y.$ In this case
we consider only partitions of the form $(n-i, i), i = 0,\dots,n$
and denote the corresponding submodule of $V^n$ by $V(n-i,i).$
This is the span of all monomials $a$ of degree $n-i$ in $x$ and
$i$ in $y.$  An $a \in V(n-i,i)$ will be called {\it homogeneous
of weight} $|a|: = n-2i;$  one then has $Ha = (n-2i)a = |a|a.$
It is easy to see (and will in any case be shown in the quantized
case) that $\ker X| V(n-i,i) \ne 0$ if and only if $n-2i \ge 0$
in which case we denote this kernel by $V(n-i,i;0).$  The
distinct simple $sl_2$ modules are the symmetric powers $V^{\odot
r}, r = 0,1, \dots$ of $V$ defined as follows:   For $r = 0$ this
is the ``trivial'' one-dimensional module annihilated by all
elements of $sl_2,$ for $r = 1$ it is $V$ itself, the ``vector
representation'', and for larger $r$ it is the $r+1$--dimensional
space spanned by the ordinary monomials (i.e., in commuting
variables) $x^{r-i}y^i, i = 0,\dots,r$ on which $H, X, Y$ act as
$x\partial_x - y\partial_y, x\partial_y, y\partial_x,$
respectively. The eigenvalue of the Casimir on $V^{\odot r}$ is
$\frac 12(r+1)^2.$  These are different for non-negative values
of $r,$ so the decomposition of any finite-dimensional $sl_2$
module into eigenspaces of the Casimir is already a decomposition
into isotypical $sl_2$ submodules.

The various $V^{\odot r}$ can be distinguished by the index of
nilpotence of $Y$ (or $X$) when it acts on them, the index being
$r+1.$  The index of nilpotence of $Y$ on $V^n$ is $n+1,$ so the
highest $V^{\odot r}$ is that with $r = n;$ it obviously occurs
exactly once and is generated by its highest weight element,
$x^n.$  Note that $V(n-i,i;0),$ the kernel of $X$ in $V(n-i,i),$
consists of all the highest weight vectors for those submodules
of $V^n$ which are isomorphic to
$V^{\odot (n-2i)}.$ Setting $V(n-i-r,i+r;r) := Y^rV(n-i,i;0)$ it
follows that $Y\,|\,V(n-i-r+1,i+r-1;r-1) \to V(n-i-r,i+r;r)$ is
an isomorphism for $r=0,\dots,n-2i-1$ and is the zero mapping for
$r=n-i.$  Therefore, every simple $S_n$--submodule of $V^n$ is
isomorphic to a submodule of some $V(n-i,i;0),$ but the
latter, we will see, are all simple. The essential step in
proving both this and that $\Bbb QS_n$ operates on each simple
module as its full ring of linear endomorphisms is to see how
$V(n-i,i;0)$ is constructed from simple $S_{n-1}$ submodules of
$V^{n-1}.$  Observe that if $n-2i > 0,$ so $V(n-i-1,i;0)$ is
defined, then $a \in V(n-i-1,i;0)$ implies $xa \in V(n-i,i;0).$
Now suppose that $b \in V(n-i,i-1;0).$ Since by hypothesis $Xb =
0$ we have $XYb = Hb = |b|b,$ where $|b| = n-2i+1 > 0.$ It
follows that $X(yb-|b|^{-1}xYb) =0,$ so we make the following
definition, generalized later.
$$\Cal Pb = yb-|b|^{-1}xYb.$$
Then we have also $\Cal Pb \in V(n-i,i;0).$  View $S_{n-1}$ as
the subgroup of $S_n$ permuting $2,\dots,n$ and leaving $1$
fixed. Then it is clear that left multiplication by $x$ viewed as
a mapping $V(n-i-1,i;0) \to V(n-i,i;0)$ (where $n-2i >0$) and
$\Cal P: V(n-i,i-1;0) \to V(n-i,i;0)$ (where $n-2i \ge 0$) are
both $S_{n-1}$--module monomorphisms.

\proclaim{Theorem 2.5} Suppose that $n-2i \ge 0.$ Then $V(n-
i,i;0)$ is the orthogonal direct sum of $\Cal PV(n-i,i-1;0)$ and
$xV(n-i-1,i;0),$ the latter being omitted if $n-2i = 0.$ Moreover
$V(n-i,i;0)$ is simple and $\Bbb QS_n$ operates on it as its full
ring of linear endomorphisms. \endproclaim

\demo{Proof} First we must show that every element $c \in V(n-
i,i;0)$ actually has the form $c = xa +\Cal Pb$ for some $a \in
V(n-i-1,i;0)$ and $b \in V(n-i, i-1;0).$  Write $c = xc_0 +
yc_1.$  Since $Xc = xXc_0 +xc_1 +yXc_1 = 0$ we must have $Xc_1 =
0,$ so $c_1 \in V(n-i,i-1;0).$  Then $c-\Cal Pc_1 =
x(c_0+|c_1|^{-1}Yc_1), |c_1| = n-2i+1,$ so the desired
$a=c_0+|c_1|^{-1}Yc_1.$  Now make the inductive assumption that
the assertion is true for all smaller values of $n,$ there being
nothing to prove when $n=1.$  Suppose that $a \in V(n-i-1,i;0)$
and $b \in V(n-i, i-1;0.$ Clearly $xa$ is orthogonal to $yb$ so
it will be orthogonal to $\Cal Pb$ if the inner product $(xa,
xYb)$ vanishes.  But $(xa, xYb) = (a, Yb)$ and $a \in V(n-
i-1,i;0)$ while $Yb \in YV(n-i,i-1;0) = V(n-i-1,i;1).$ There are
two ways to see that these are orthogonal. First, they lie in
different eigenspaces of the Casimir.  Without invoking the
Casimir, however, observe that by the inductive hypothesis $xV(n-
i-1,i;0)$ and $\Cal PV(n-i,i-1;0)$ are non-isomorphic simple
$\Bbb QS_{n-1}$--submodules of $V(n-i,i;0)$ considered as an
$S_{n-1}$--module. They are therefore orthogonal since, as
remarked, the orthogonal projection of one module on another is a
submodule of the second.  To see that $\Bbb QS_{n-1}$ acts on
each as its full ring of linear endomorphisms,  suppose that
$V(n-i-1,i;0)$ has rank $r$ and $V(n-i,i-1;0)$ has rank $s.$
These are the same as the ranks of their images in $V(n-i,i;0),$
so taking a basis of $V(n-i,i;0)$ formed by combining bases of
these images, one can view the representation of $\Bbb QS_n$ on
$V(n-i,i;0)$ as one by $(r+s)\times (r+s)$ matrices which already
contains the direct sum of the $r \times r$ and $s \times s$
matrices. But neither  $xV(n-i-1,1;0)$ nor $\Cal PV(n-i,i-1;0)$
is an $S_n$ submodule. It follows that we must have the full ring
of $(r+s) \times(r+s)$ matrices. $\square$ \enddemo

The orthogonality of the summands out of which $V(n-i,i;0)$ is
built gives the inductive construction of its canonical
orthogonal basis: combine the images of those of $V(n-i,i-1;0)$
and (if $n-2i > 0$) $V(n-i-1,i;0)$ under $\Cal P$ and left
multiplication by $x,$ respectively. The elements of the
resulting basis may therefore be indexed by sequences of length
$n$ in $x$ and $\Cal P,$ where $x$ appears $n-i$ times and $\Cal
P$ appears $i$ times and where, in each terminal segment
(subsequence consisting of the last $j$ elements for each $j$)
the number of $\Cal P$'s does not exceed the number of $x$'s.
Instead of sequences we may take the index set  to be two-rowed
Young diagrams.  For if the sequence (which must end with $x$) is
given then the associated Young diagram is constructed as
follows: Put  ``1'' in the first row and first column. The
position of ``2'' is determined by the next--to--last symbol in
the sequence; if it is an $x$ then put ``2'' in the next position
in the first row, if it is a $\Cal P$ put it in the first
position in the second row.  If the integers through $i$ have
been positioned using the last $i$ entries in the sequence, then
$i+1$ is put in the first open position of the first or second
row according as the $n-i$th entry in the sequence is an $x$ or a
$\Cal P.$  The condition on the sequence is precisely that this
should give a Young diagram, so we have

\proclaim{Proposition 2.6} The canonical orthogonal basis of
$V(n-i,i;0)$ is indexed by the two-rowed Young diagrams in which
the first row has length  $n-i$ and the second row has length
$i.$ In particular, the rank of the module is the number of such
diagrams. $\square$ \endproclaim

In view of this, we will call a sequence of $x$'s and $\Cal P$'s
in which every terminal segment contains no more $\Cal P$'s than
$x$'s a {\it generating sequence.}  For two variables, the number
of generating sequences is the difference of the binomial
coefficients $\binom ni - \binom n{i-1},$ which we denote by
$\langle n,i\rangle.$ The recursion formula for these is
essentially the same as for the binomial coefficients themselves:
$\langle n,i\rangle = \langle n-1,i-1\rangle + \langle
n-1,i\rangle,$ where $n-2i \ge 0$ and the second summand is
omitted if $n-2i = 0.$

As an elementary example, consider the simple two-dimensional
representation  $V(2,1;0)$ of $S_3.$ Its orthogonal basis is
indexed by the two sequences $x\Cal Px, \Cal Pxx.$ The first
basis element is thus $x(\Cal Px) = x(yx-xy)$ and the second is
$\Cal P(xx) = yxx-\frac 12 x(yx+xy).$ The two basis elements
correspond, respectively, to the Young diagrams having $1,2$ in
the first row (second row empty) and having $1,2$ in the first
column (second column empty).  It is important, as we shall prove
later in the quantized case, that the norm of any element of the
canonical orthogonal  basis of any $V(n-i,i;0)$ (i.e., the sum of
the squares of its coefficients when expressed in terms of the
original monomial basis) must be a unit in the ring $\Bbb
Z[1/n!].$ In the small example with $n=3$ just given, the norms
of the two basis elements are 2 and 3/2, respectively.
For the case of arbitrary $d$ we shall have to generalize $\Cal
P$ to a sequence of $d$ operators. (Here we actually have two;
the first is left multiplication by $x.$)  We turn now to the
quantized case.

\subhead 3. Representation of the Hecke algebra on $V^n$
\endsubhead
Recall that the Hecke algebra $\Cal H_n$ of $S_n,$ whose
representations we must study, is the algebra generated over
$\Bbb Z[q, q^{-1}]$ by elements $T_{s_i}, i = 1, \dots, n-1$
corresponding to the generators $s_i = (i,i+1)$ of $S_n.$ Its
multiplication is given by $T_sT_w = T_{sw}$ if $s$ is one of
these generators and $w \in S_n$ is an element with length
$\ell(sw) > \ell(w),$ and $T_s^2 = (q-q^{-1})T_s +1.$ It is a
free module of rank $n$ over $\Bbb Z[q,q^{-1}].$  The familiar
``Artin presentation'' of $S_n$ is by generators  $s_i, i =
1,\dots, n-1$ with
$$\align
s_is_{i+1}s_i &= s_{i+1}s_is_{i+1}\\
s_is_j &= s_js_i \qquad \text{if} \qquad |i-j| > 1\\
s_i^2 &= 1
\endalign$$
It follows that for the special case of $S_n$ we obtain an
equivalent definition of the Hecke algebra $\Cal H_n$ by
requiring  that the generators $T_i$ satisfy the ``braid'' and
``Hecke''  relations
$$\align
T_iT_{i+1}T_i &= T_{i+1}T_iT_{i+1}\\
T_iT_j &= T_jT_i \qquad \text{if} \qquad |i-j| > 1\\
T_i^2 &= (q-q^{-1})T_i + 1
\endalign$$
For the Hecke algebra of $S_n$ certainly satisfies these
relations so if the coefficient ring were a field it would be
sufficient to show that the algebra these define has no greater
dimension than that of the Hecke algebra. For this it is
sufficient to  show that if a formal product $w$ of generators of
$S_n$ gives a non-reduced expression for some element of $S_n$
then the corresponding product $\tau$ of generators $T_i :=
T_{s_i}$ satisfying the second set of relations can also be
shortened.  But $w$ can be shortened only if by applying the
braid relations it can be rewritten to contain the square of a
generator, in which case $\tau$ can be shortened also.  Over a
domain we must note further only that the Hecke algebra, which is
in any case a quotient of the algebra defined by the above
relations, is a free module, so the quotient map splits.  Many
authors, e.g., [W] simply adopt the second definition.

Now let $V$ be a free module of rank $d$ over $\Bbb Z[q,q^{-1}].$
There is then a natural representation of $\Cal H_n$ on $V^n$
given as follows. The standard $d^2 \times d^2$ quantum Yang--
Baxter matrix
$$R =\sum \Sb i \ne j\\ i,j=1\endSb^d e_{ii}\otimes e_{jj} +
q\sum_{i=1}^d e_{ii}\otimes e_{ii} + (q-q^{-1}) \sum_{1\le j <
i \le n} e_{ij}\otimes e_{ji} $$
for the simple Lie algebra $sl_d$ (cf. \cite{FRT}) may be viewed
as having coefficients in $\Bbb Z[q,q^{-1}].$ We view this as
operating on $V^2$ so: Let the basis of $V$ be $x_1,\dots,x_d,$
so that of $V^2$ consists of the $x_ix_j (= x_i\otimes x_j)$ in
lexicographic order. Then set $(e_{ij}\otimes e_{kl})x_rx_s =
x_ix_k$ if $j=r, l=s$ and 0 otherwise.  Let $(12)$ operate as the
interchange of tensor factors in $V^2.$ As a matrix, we have  $$
(12) = \sum_{i,j}e_{ij}\otimes e_{ji}.$$
Set
$$\overline R = (12)R
   =\sum \Sb i \ne j\\ i,j=1\endSb^d e_{ji}\otimes e_{ij} +
q\sum_{i=1}^d e_{ii}\otimes e_{ii} + (q-q^{-1}) \sum_{1\le j <
i \le n} e_{jj}\otimes e_{ii}.$$
This is symmetric. (The transpose of $e_{ij}\otimes e_{kl}$ is
$e_{ji}\otimes e_{lk},$ the pair $(i,k)$ being the row index and
$(j,l)$ the column index.)  It is also ``{\it balanced}\,'',
i.e., in every term the sum of the column indices equals the sum
of the row indices and it is therefore a direct sum of matrices
in each of which these ``{\it weights\,}'' are constant.  (This
concept is meaningful for an arbitrary tensor power of a matric
algebra; for the first power, a balanced matrix is simply
diagonal.)

Let $\overline R_i, i =1,\dots,n-1$ denote the operation of
$\overline R$ in tensor factors $i, i+1$ of $V^n$.  It is a basic
fact that these satisfy the braid relations, cf. \cite{FRT}:
$$\align \overline R_i\overline R_{i+1}\overline R_i
 &= \overline R_{i+1}\overline R_i\overline R_{i+1}\\
\overline R_i\overline R_j &= \overline R_j\overline R_i \qquad
\text{if} \qquad |i-j| >1.\endalign$$
Sending $T_i$ to $\overline R_i$ therefore induces a
representation of $\Cal H_n$ on $V^n.$  Note that $\overline R$
is self-adjoint -- its matrix is symmetric relative to the
standard basis of $V^2$ consisting of the $x_ix_j$ in
lexicographic order --  and therefore so are all the $\overline
R_i.$  Here, for example, is the operation of $\overline R$ on
$V^2$ in the case where $d=2$:
$$\align
\overline Rx^2 &= qx^2 \\
\overline Rxy &= (q-q^{-1})xy + yx\\
\overline Ryx &= xy\\
\overline Ry^2 &= qy^2
\endalign$$
As in the classical case, we have the following basic
\proclaim{Lemma 3.1} The orthogonal projection of one $\Cal
H_n$--submodule of $V^n$ on another is a submodule of the
second.\endproclaim
\demo{Proof} This would be trivial if, as in the classical case,
the representation were generated by orthogonal transformations.
However, viewing $q$ as a real parameter there are constants $C$
and $S$ such that $C\overline R + S$ is orthogonal. In fact, set
$$\gamma = \sum_{i<j}e_{ij}\wedge e_{ji} \qquad
\text{where}\qquad
e_{ij}\wedge e_{kl} = \frac12(e_{ij}\otimes e_{kl}-e_{kl}\otimes
e_{ij}).$$
This is the infinitesimal of the deformation from
$Usl_n$ to $U_qsl_n.$ Set $q = \sec t - \tan t.$
Then $(\cos t)\overline R+\sin t=e^{-t\gamma}(12)e^{t\gamma}$ in
which all the factors are orthogonal; here $C = \cos t, S = \sin
t.$ (Cf. \cite{GGSk2} where the interchange of tensor factors
``$(12)$''is denoted by $P.$) It follows in this case that the
representation of the Hecke algebra is generated by orthogonal
transformations of $V^n,$ so here the orthogonal projection of
one $\Cal H_n$--submodule on another is indeed a submodule of the
second.  This, however, is a purely formal property and therefore
holds generally. $\square$\enddemo

The Hecke algebra may be viewed as a ``quantization'' of the
group ring $\Bbb ZS_n$ and the commutant of its representation on
$V^n$ is a representation of a particular quantization
$U_q:=U_qsl_d$, so we describe that which we use starting with
the case $d = 2.$  Observe first that the unquantized $U$ is a
Hopf algebra with primitive generators $H,X, $ and $Y.$ The
Drinfel'd-Jimbo quantization (cf. \cite{D, J1, J3}) replaces $H$
by a pair of invertible group--like generators $K,K^{-1}$ (i.e.,
$\Delta K = K\otimes K, \Delta K^{1} = K^{-1} \otimes K^{-1}).$
We want, in effect, that $K =q^H,$  so $K$ will act on the tensor
algebra $TV$ as the automorphism sending $x$ to $qx$ and $y$ to
$q^{-1}y,$ which forces the commutation relations  $KX = q^2XK,
KY = q^{-2}YK.$  The remaining multiplication and
comultiplication rules are given by
$$\gather
qXY-q^{-1}YX = (q^{-1}-q)^{-1}(1-K^2) \\
\Delta X = X\otimes1 + K\otimes X,\qquad
\Delta Y =  Y\otimes1 + K\otimes Y
\endgather$$
(For a full definition one needs in addition, the $q-$analogues
of the Serre relations giving, in particular, the nilpotence of
$\operatorname{ad} X$ and $\operatorname{ad} Y$, but we do not
need these here, nor shall we need the quantized antipode.)  The
right side of the first relation, which has the correct
quasi--classical limit
 as $q \to 1$, namely $H,$ may seem to pose a
problem since we do not assume that $q-q^{-1}$ is invertible. We
therefore take the left side, which is well-defined, as an
additional generator denoted simply by $H.$  This ``quantized''
$H$ should be distinguished from the original unquantized one
which, if we should need it, would be denoted $H_0.$
Both $K$ and the quantized $H$ will be seen to be well-defined
operators on $V^n.$

The operations of $X, Y$ on $V$ itself are the same as before
quantization, that is, $Xx= 0, Xy = x, Yx = y, Yy = 0$. However,
the new comultiplication (which defines the tensor product of
modules, in particular extending the operation of $U_q$ to $V^n$
for every $n$) now specifies that for homogeneous elements
$\alpha,\beta \in V^n$ we have
$$X(\alpha\beta) = (X\alpha)\beta+q^{|\alpha|}\alpha X\beta \quad
Y(\alpha\beta) = (Y\alpha)\beta+q^{|\alpha|}\alpha X\beta.$$ From
the  definitions one can also readily deduce that
$$H\alpha = q|\alpha|_{q^2}\alpha$$
It follows, in particular, that if $X\alpha = 0$ then $XY\alpha
=|\alpha|_{q^2} \alpha.$   Most important, as the reader should
check, the operations of $U_q$ so defined on $V^n$ commute with
those of the $\overline R_i$ and therefore are $\Cal H$--module
morphisms.
The second cohomology of a simple Lie algebra taken with
coefficients in itself vanishes, so a simple Lie algebra admits
only trivial deformations of either itself as a Lie algebra or of
its universal enveloping algebra, giving the following quantized
form of Schur's theorem in the generic case, cf. \cite{J2}.
\proclaim{Theorem 3.2} Setting $q=1+t$ and extending coefficients
to $\Bbb Q[[t]],$ the operations of $\Cal H_n$ and $U_q$ on $V^n$
are mutual commutants. \endproclaim
\demo{Proof} This is an exercise in deformation theory: Note that
after the extension $\Cal H_n$ is isomorphic to $\Bbb QS_n[[t]]$
and, because $sl_2$ is simple, $U_q$ likewise becomes isomorphic
to $U$ with coefficients extended. $\quad\square$ \enddemo

Since $U_qsl_d$ is a trivial deformation of $Usl_d$ it has (after
extension of coefficients) a Casimir operator. Although we shall
not need it, for $d=2$ and our specific quantization, this is
given by
$$\Cal C = (q-q^{-1})^{-2}(qK+q^{-1}K^{-1}-2)+q^{-1}K^{-1}YX.$$
(This is adapted from Rosso, \cite{R}.)  The strange form is
forced by our quantization, which in turn is forced by the
requirement that $U_q$ be in the commutant of $\Cal H.$ It is
easy to verify directly that $\Cal C$ is in the center of $U_q.$
Its limit as $q \to 1$ is half our previous classical Casimir
operator, which we now would denote by $\Cal C_0.$ Recall that we
still have a symmetric inner product on $V^n.$

\proclaim{Theorem 3.3} The quadratic and higher Casimir operators
are self-adjoint.
\endproclaim
\demo{Proof} The assertion will hold if it does so generically,
and is not affected by extension of coefficients, so we may
replace $q$ by $1+t$ and take coefficients to be $\Bbb Q[[t]].$
The Casimir is then a central element of $\Cal H_n$ which, since
coefficients have been extended, is isomorphic to $\Bbb
QS_n[[t]].$ The generators of $\Cal H_n$ were self-adjoint.
Therefore transposition carries it as a whole into itself. Since
the center is preserved, transposition can only act as a
permutation on the central idempotents.  But when $t=0$ the
central elements are self-adjoint, so this permutation, which is
a continuous function of $t$, is the identity when $t=0.$
Therefore it is the identity on the center for all $t,$ so the
Casimir is preserved. $\square$ \enddemo

The case $d=2$ proceeds exactly as in the classical case. Recall
that $\Bbb Z_{q,n}:=\Bbb Z[q,q^{-1},1/n_{q^2}!]$ which we
henceforth always take as coefficients. The quantized $H$
operates semisimply on $V^n$ which is the direct sum of its
eigenspaces $V(n-i,i);$ the eigenvalue on this is
$q(n-2i)_{q^2},$ which is a unit. In fact, any eigenvalue of $H$
on a module of finite rank must be of the form $q\lambda_{q^2}$
with $\lambda$ a non-negative integer. For if $\alpha$ is an
eigenvector and we write the eigenvalue formally as
$q\lambda_{q^2}$ with
$\lambda_{q^2} = (1-q^{2\lambda)}/(1-q^2),$ then it is easy to
verify that $HX^m\alpha = q(\lambda+2m)_{q^2}\alpha$ and
$HY^m\alpha = q(\lambda-2m)_{q^2}\alpha.$ A simple telescoping
induction on the relation \linebreak $qXY - q^{-1}YX = H$ gives
$$q^mXY^m - q^{-m}Y^mX = q^{m-1}HY^{m-1} +
q^{m-3}YHY^{m-2} +  \dots +q^{-m+1}Y^{m-1}H,$$
so we have
\proclaim{Lemma 3.4} If $H\alpha = q\lambda_{q^2}\alpha$     and
$X\alpha
= 0$ then\newline \phantom{............}\qquad $  XY^m\alpha =
q^{-2m+2}(\lambda-m+1)_{q^2}m_{q^2}Y^{m-1}\alpha. \qquad\square$
\endproclaim

The coefficient on the right will not vanish for any $m$ unless
$\lambda$ is a non-negative integer, which therefore must be the
case if the module is of finite rank.  Moreover, one then has
that if $X\alpha = 0$ and $H\alpha = 0$ then $Y\alpha = 0$ as
well.  The simple modules of finite rank over $U_q$ thus look
precisely like those over $U$, provided the non-zero
coefficients which appear in the analogues of the usual formulas
are units, which is exactly what we have supposed.

Returning to $V^n$, the submodule $\ker X$  is obviously
homogeneous, i.e., the direct sum of its components in each of
the $V(n-i,i)$. As before, write $\ker X|V(n-i,i) = V(n-i,i;0)$.
If $\alpha \in V(n-i-1,i;0)$ then $x\alpha \in V(n-i,i;0)$; if
$\beta \in V(n-i,i-1;0)$ set $\Cal P\beta = -
(q|\beta|_{q^2})^{-1}xY\beta + y\beta$ for $|\beta| \neq 0$
and  $\Cal P\beta = 0$ if $|\beta| = 0.$  (This is the quantized
version of our earlier $\Cal P.$) It is easy to check that this
is also in $V(n-i,i;0).$  Exactly as in the unquantized case, set
$V(n-i,i;r) = Y^rV(n-i+r,i-r;0) $ for $r = 0,\dots, i$ (they
will be seen to vanish for larger $r$).  These are all $\Cal H_n$
submodules.

\proclaim{Lemma 3.5}
\roster
\runinitem If $0\le r\le n-2i$ then
$X:V(n-i-r,i+r;r) \to V(n-i-r-1,i+r+1;r+1)$ and
$Y:V(n-i-r-1,i+r+1;r+1) \to V(n-i-r,i+r;r)$ are isomorphisms and
the compositions $XY$ and $YX$ are multiplication by an
invertible constant.
\item $Y^rV(n-i,i;0) = 0$ for $r > n-2i$, so $V(n-i-r,i+r;r)=0$
for $r>n-2i,$ and $V(n-i,i;0)=0$ for $n-i<i.$
\endroster
\endproclaim
\demo{Proof}
(1) Suppose that $1 \le r \le n-2i$ and that $\alpha \in V(n-
i,i;0),$  so $\beta=Y^{r-1}\alpha \in V(n-i-r+1,i+r-1;r-1).$
Since $X\alpha=0,$ by Lemma 3.4 we have
$$XY\beta=XY^r\alpha=q^{-2r+2}(|\alpha|-r+1)_{q^2}r_{q^2}\beta.$$
 Then the coefficient  $c:=q^{-2r+2}(|\alpha|-r+1)_{q^2}r_{q^2}$
is  invertible, so  $XY$ is an automorphism, whence $X$ is onto
and $Y$ is one--to--one. However, by definition $Y$ is onto, thus
invertible, so $X$ and $Y$ are both isomorphisms. Now $XY\beta =
c\beta$ implies $YXY\beta = cY\beta.$ Since $Y$ is onto, $YX$ is
just multiplication by the same $c.$

(2) Clearly some power of $Y$ annihilates $V(n-i,i;0).$ If $r$ is
the first such then Lemma 3.4 with $\lambda = n-2i$ implies that
$n-2i-r+1 = 0,$ else $Y^{r-1}$ would already annihilate $V(n-
i,i;0).\quad\square$
\enddemo

Using the quantized $\Cal P$ and exactly the same proof as in the
classical case, we now have

\proclaim{Lemma 3.6} If $n-2i \ge 0$ then
 $V(n-i,i;0)$ is the orthogonal direct sum of $xV(n-i-1,i;0)$ and
$\Cal P V(n-i,i-1;0),$ where the first summand appears only if
$n-2i > 0. \quad\square$
\endproclaim

Again, with arguments identical to those in the classical case we
have, when $d=2,$

\proclaim{Theorem 3.7} Every simple $\Cal H_n$--submodule of
$V^n$ is isomorphic to one of the $V(n-i,i;0).$ Each of the
latter has a canonical orthogonal basis indexed by the two-rowed
Young diagrams with rows of length $n-i$ and $i,$ and $\Cal H_n$
acts on each as its full ring of linear endomorphisms. $\square$
\endproclaim

It follows, as before, that the norms of the basis elements are
units. Also, $V(n-i,i)$ is the orthogonal direct sum of its
submodules $V(n-i,i;r),$ each of which is simple.

\subhead 4. The case of arbitrary rank
\endsubhead

Suppose now that $V$ is a free module of rank $d$ over $Z_{q,n}=
\Bbb Z[q,q^{-1}, 1/n_{q^2}!]$ with basis $x_1,\dots,x_d$.  The
``standard'' quantization $U_q = U_qsl_n$ (cf. \cite{D}) may be
viewed as generated by elements $X_i, Y_i, K_i, H_i, i =
1,\dots,d-1$ which for fixed $i$ generate a subbialgebra
isomorphic to $U_qsl_2.$  Here $X_i$ and $Y_i$ act on $V$ like
$x_i\partial/\partial x_{i+1}$ and $x_{i+1}\partial/\partial
x_i$, respectively. These act on $V^n$ so: If  $\alpha$ is
homogeneous in variables $x_i$ and $x_{i+1}$ of degrees $n_i$ and
$n_{i+1}$ in each, respectively, then  set
$|\alpha|_i = n_i - n_{i+1}$  and set
$$X_i(\alpha\beta) = (X_i\alpha)\beta +
q^{|\alpha|_i}\alpha(X_i\beta); \quad Y_i(\alpha\beta) =
(Y_i\alpha)\beta + q^{|\alpha|_i}\alpha(Y_i\beta).$$
In particular, $X_i$ and $Y_i$ treat all variables $x_j$ with
$j  \ne i,i+1$ as constants.  It is then evident from the
preceding formula that the $X_i$ and $Y_i$ act as $\Cal H_n$
morphisms under the induced representation of $\Cal H_n$ on
$V^n$.  Moreover, if  $|i-j| > 1$ then $X_i$ and $Y_j$ commute,
as do $X_i$ and $X_j$, and $Y_i$ and $Y_j,$ but unlike the
classical case, $X_i$ and $Y_{i\pm1}$ only ``$q$--commute'',
i.e., we have $$X_iY_{i\pm1} = qY_{i\pm1}X_i.$$
Each $X_i, Y_i, K_i$ and (quantized) $H_i$ together generate a
subalgebra of $U_qsl_n$ isomorphic to the subalgebra of $U_qsl_2$
generated by $X, Y, K$ and $H.$

If $\frak p = (n_1,\dots,n_d)$ is a non-increasing partition of
$n$ into $d$ parts (some of which may be zero) then we will write
$\frak p! = n_1!n_2!\cdots n_d!$ and the usual multinomial
coefficient $n!/\frak p!$ will be denoted simply by $\binom
n{\frak p}$.  When $\frak p'$ is a non-increasing   partition of
$n-1$ into $d$ parts differing in exactly one place from $\frak
p$ then we write $\frak p' \to \frak p$. Denote by $V(\frak p)$
the submodule of $V^n$ spanned by all monomials of degree $n_i$
in $x_i, i = 1,\dots,d$, where   the $n_i$ are the parts of
$\frak p$.  (This is meaningful even if $\frak p$  fails to be
non--increasing.)  In analogy with the preceding section, we set
$$ V(\frak p; 0) = \bigcap_{i = 1}^{d-1} \ker X_i.$$
(This, too, is meaningful even if $\frak p$ fails to be non--increasing,
but from the preceding section it then vanishes.)
For the rest of this section, unless specified, we assume the
coefficient ring to be $\Bbb Z_{q,n} = \Bbb Z[q,q^{-1},
1/n_{q^2}!].$ We can now construct the $V(\frak p;0)$ from the
$V(\frak p';0)$ with $\frak p' \to \frak p$ in a way generalizing
that of the preceding section. Suppose that $\frak p =
(n_1,\dots,n_d)$ and that $\frak p' =
(n_1,\dots,n_r-1,\dots,n_d)$ where the latter is still non-
increasing.  Define inductively polynomials in
$Y_1,\dots,Y_{r-1}$ (whose dependence on $\frak p$ and $r$ we
momentarily suppress) by $\widehat P_0 = 1$, $\widehat P_1 =
Y_{r-1}$, and
$$\widehat P_i = (n_{r-i+1}-n_r+i)_{q^2}Y_{r-i}\widehat
P_{i-1}-q(n_{r-i+1}-n_r +i-1)_{q^2}\widehat P_{i-1}Y_{r-i}, \quad
1 < i < r.$$
It follows that $\widehat P_i$ is a polynomial in the  operators
$Y_{r-i}, Y_{r-i+1}, \dots , Y_{r-1}$ only, and is  homogeneous
of degree 1 in each.  Applied to a homogeneous $\alpha$, it
decreases its degree in $x_{r-i}$ by one and increases its degree
in $x_r$ by one.  Set $P_0 = 1$ and for $ i > 0$ set
$$P_i = (-1)^iq^{-i}[(n_{r-1}-n_r+1)_{q^2}\dots(n_{r-i}-n_r
+i)_{q^2}]^{-1}\widehat P_i.$$
Note that the ``arguments'' of the coefficients never exceeds
$n$.  Finally, set
$$\Cal P(\frak p' \to\frak p) = x_rP_0 + x_{r-1}P_1 + \dots +
x_1P_{r-1}$$
which we view as a mapping $V(\frak p') \to V(\frak p)$.
\proclaim{Lemma 4.1}
Suppose $\frak p'= (n_1,\dots,n_r-1,\dots,n_d) \to \frak p =
(n_1,\dots,n_r,\dots,n_d).$
\roster
\item For all $i = 1,\dots,r-1,$ if $\alpha\in V(\frak p')$ is in
$\ker X_j,$ then
$$\align
X_j\widehat P_i\alpha &= 0\qquad\text{for}\qquad j\neq
r-i,\tag a\\
X_{r-i}\widehat P_i\alpha &= (n_{r-i} - n_r +i)_{q^2}\widehat
P_{i-1}\alpha. \tag b
\endalign$$
\item If $\alpha \in V(\frak p';0) = \bigcap_{j=1}^{d-1} \ker
X_j\bigm| V(\frak p')$ then $\Cal P(\frak p \to\frak p')\alpha
\in V(\frak p;0).$
\item $V(\frak p;0) = \sum_{\frak p' \to \frak p} \Cal P(\frak
p'\to \frak p)V(\frak p';0).$
\endroster
\endproclaim
\demo{Proof}
(1) When $i=1$ statement (a) is immediate because $X_j$ commutes
with  $Y_{r-1}$ and $X_j\alpha=0;$ (b) asserts simply that if
$X_{r-1}\alpha = 0$ then $X_{r-1}Y_{r-1}\alpha =
(n_{r-1} - n_{r}+1)_{q^2}\alpha.$ This holds from the
corresponding  assertion for $d = 2.$ Suppose now that $i > 1$
and that the  assertions hold for all smaller $i.$ For
simplicity, write  $(n_{r- i+1} - n_r +i-1) = \lambda$ . Then
$$X_j\widehat P_i\alpha  = (\lambda + 1)_{q^2}X_jY_{r-i}\widehat
P_{i-1}\alpha - q\lambda_{q^2}X_j\widehat P_{i-1}Y_{r-i}\alpha.
\tag*$$
Suppose first that $j > r-i+1$ or $j<r-i-1.$ Then $X_j$ commutes
with  $Y_{r-i}$ so the first term vanishes by the inductive
hypothesis.  So does the second, since $Y_{r-i}\alpha$ is still
in $\ker X_j.$ For $j = r-i+1$, using the commutation relations
and  the inductive hypothesis, including that on (b), one finds
that  the first term is
$$q(\lambda+1)_{q^2} Y_{r-i}X_{r-i+1}\widehat P_{i-1}\alpha =
q(\lambda+1)_{q^2}\lambda_{q^2}Y_{r-i}\widehat P_{i-2}\alpha.$$
Since $Y_{r-i}$ commutes with all the $Y$ in $\widehat P_{i-2}$,
to show that the whole expression vanishes we must show only that
$X_{r-i+1}\widehat P_{i-1}Y_{r-i}\alpha =
(\lambda+1)_{q^2}\widehat P_{i-2}Y_{r-i}\alpha.$
As before, $Y_{r-i}\alpha$ is still in $\ker X_{r-i+1},$ so we
may again apply the inductive hypothesis on  (b) noting, however,
that the degree of $Y_{r-i}\alpha$ in $x_{r-i+1}$ is greater by 1
than the corresponding degree of $\alpha$. When $j=r-i-1,$ \
$X_j$ commutes or $q$--commutes with $\widehat  P_{i-1}$ and
$Y_{r-i},$ so $X_j\alpha=0$ implies that in (*) both factors are
zero.  This proves (1a). For (b), setting $j = r-i$ in (*), note
that $X_{r-i}$ annihilates  $\widehat P_{i-1}\alpha,$ the degree
of which in $x_{r-i+1}$ is one less than that of $\alpha.$    If
we let $\mu=n_{r-i}-n_{r-i+1}$ be the $r-i$ norm of $\alpha,$
then
$$\gather
X_{r-i}Y_{r-i}\alpha=q^{-1}H_{r-i}\alpha=\mu_{q^2}\alpha\qquad
\text{and}\\
X_{r-i}Y_{r-i}\widehat P_{i-1}\alpha=q^{-1}H_{r-i}\widehat
P_{i-1}\alpha=(\mu+1)_{q^2}\widehat P_{i-1}\alpha.
\endgather$$
The first term on the right of (*)is
$(\lambda+1)_{q^2}(\mu+1)_{q^2} \widehat P_{i-1}\alpha.$ For
the second term,   since $\widehat P_{i-1}$ is homogeneous of
degree 1 in $Y_{r-i+1}$ and all its other factors commute with
$X_{r-i}$ one gets
$-q^2\lambda_{q^2}\mu_{q^2}\widehat P_{i-1}\alpha.$  That the sum
of the coefficients of $\widehat P_{i-1}\alpha$ which appear is
indeed $(\lambda+\mu+1)_{q^2}$ is an exercise in the definition
of the ``q--numbers". (Note that it is true when $q=1.$)

(2) Suppose that $\alpha \in V(\frak p'; 0)$ and consider
$X_jx_iP_{r-i}\alpha.$  If $j \neq i, i+1$ then this is a
multiple of $x_iX_jP_{r-i}\alpha$ and therefore vanishes, by
(1a).  But $X_i(x_{i+1}P_{r-i-1} + x_iP_{r-i})\alpha =
(x_iP_{r-i-1} +q^{-1}x_{i+1}X_iP_{r-i-1} +qx_iX_iP_{r-i})\alpha.$
The middle term vanishes by (1a), so this is $x_i(P_{r-i-1} +
qX_iP_{r-i})\alpha$ which vanishes from (1b) and the definition
of the $P_i.$

(3) Every  $b \in V(\frak p; 0)$ can be written in the form
$x_1a_1 + \dots +x_da_d.$  If $a_r \neq 0$ while $a_{r+1} =
\dots = a_d = 0$ then we say that $b$ has {\it length} $r.$  It
is sufficient to show that $b$ can be written as the sum of an
element of $\sum \Cal P(\frak p' \to \frak p)V(\frak p';0)$
and one of shorter length.  Now $X_ia_r = 0$ for all $i$ since
$X_ib$ contains a non-zero multiple of $x_rX_ia_r$ and no other
term beginning with $x_r$. Since $a_r\neq 0$ yet lies in $\cap
\ker X_i,$  it follows that $n_r > n_{r+1}$ so $\frak p' =
(n_1,\dots,n_{r}-1,\dots,n_d) \to \frak p.$ But then
$b-\Cal P(\frak p' \to \frak p)a_r$ is a shorter element of
$V(\frak p; 0).\quad\square$
\enddemo

There will be no ambiguity if we now denote
$\Cal P(\frak p' \to \frak p):V(\frak p';0) \to V(\frak p;0)$ by
$\Cal P_r$ where $r$ is the unique place in which $\frak p'$ and
$\frak p$ differ. More generally, we may view it as an operator
defined on any $V(\frak p';0)$ where $\frak p'$ is a non-
increasing partition of some $n' \le n$ into $d$ parts with the
property that its $r$th part can be increased by $1$ without
losing the property of being non-increasing.  Now suppose that we
have a sequence
$$\frak p^{(n)} = (0^d) \to (1,0^{d-1}) = \frak p^{(n-1)} \to
\dots
\to \frak p''\to \frak p' \to \frak p$$
where each $\frak p^{(i)}$ is a non--increasing partition of $n-
i.$  For each ``morphism'' $\frak p^{(i)} \to \frak p^{(i-1)}$
there is a unique place  $r_i$ in which $\frak p^{(i)}$ and
$\frak p^{(i-1)}$ differ, so the sequence $(r_1,r_2,\dots,r_n)$
determines the sequence of morphisms and therefore an element
$\Cal P_{r_1}\Cal P_{r_2}\dots\Cal P_{r_n}1 \in V(\frak p;0).$
Here necessarily $r_n = 1$ so $\Cal P_{r_n}1 = x_1.$  The
sequence is not arbitrary: in any ``terminal segment''
$(r_j,r_{j+1},\dots,r_n)$ the number of times any integer $k$
appears amongst the $r_i, i = j,\dots,n,$ can not exceed the
number of times that $k-1$ appears, since the partition
$\frak p^{(j-1)}$ would otherwise fail to be non-increasing.
However, that is the only restriction.  These generating
sequences $(r) = (r_1,\dots,r_n)$ correspond to $d$--rowed Young
diagrams: the value of $r_{n-i}$ is the row in which $i+1$ is
placed, and it is put to the immediate right of any other
integers previously placed.  Since $r_n = 1,$ we must put ``$1$''
in the first position of the first row.  Then $r_{n-1},$ whose
value can be either $1$ or $2,$ determines the placement of
``$2$''. If it is $1$ then ``$2$'' goes in the first row to the
right of ``$1$'', if it is $2$ then ``$2$'' begins the second
row, etc., generalizing the discussion of the previous section.
The correspondence between generating sequences and $d$--rowed
Young diagrams is clearly a bijection.  Recall that our
coefficient ring is
$\Bbb Z_{q,n} = \Bbb Z[q,q^{-1}, 1/n_{q^2}!].$

\proclaim{Theorem 4.2} ({\it i}) Each $V(\frak p;0)$ is a simple
$\Cal H_n$ module, no two are isomorphic, and $\Cal H_n$ acts on
each as its full ring of $Z_{q,n}$--linear endomorphisms. It has
a canonical orthogonal basis indexed by the $d$--rowed Young
diagrams associated with the partition $\frak p.$ ({\it ii})
$\Cal H_n$ is a direct sum of total matric algebras. ({\it iii})
Every simple submodule of $V^n$ is isomorphic to exactly one of
the $V(\frak p;0).$ If $d \ge n$ then every simple $\Cal H_n$
module is isomorphic to one of these. \endproclaim
\demo{Proof} Make the inductive assumption that the assertions
hold for all smaller values of $n.$
(i) As an $S_{n-1}$--module, $V(\frak p;0)$ is, by the preceding,
a sum of mutually non-isomorphic simple modules which must
therefore be mutually orthogonal inside $V(\frak p;0).$ None of
these summands is an $S_n$--module and $H_{n-1}$ already acts on
each as its full ring of linear endomorphisms.  Suppose that the
ranks of these submodules are $m_1,\dots,m_s,$ so that of
$V(\frak p;0)$ is $m:=m_1+\dots+m_s.$ Then $\Cal H_n$ operates on
$V(\frak p;0)$ as a subalgebra of an $m \times m$ matric algebra
which already contains the direct sum, denote it $D,$ of $m_1
\times m_1,\dots, m_s \times m_s$ total matric algebras and in
which none of these is a submodule. Over a field it is easy to
see that $\Cal H_n$ must then operate as the full $m \times m$
matric algebra. (Pass to the skeleton of the induced operation of
$\Cal H_n$, which permits one to assume that all the $m_i = 1.$)
Hence it holds if we map $\Bbb Z_{q,n}$ into any field and take
as coefficients the subfield generated by the image, but this
implies that it holds for $\Bbb Z_{q,n}$ itself. It follows that
$V(\frak p;0)$ is a simple $\Cal H_n$ module. If any two were
isomorphic, say $V(\frak p_1;0)$ and $V(\frak p_2;0),$ then they
would already be isomorphic as $\Cal H_{n-1}$ modules. But for
distinct non-increasing partitions $\frak p$ of $n$, the sets of
$\frak p' \to \frak p$ are distinct although generally not
disjoint.  It follows that
$V(\frak p_1;0)$ and $V(\frak p_2;0)$ can not be isomorphic as
$\Cal H_{n-1}$ modules.  The assertion about the canonical
orthogonal basis is now evident from the remarks preceding the
theorem.
(ii), (iii) As in the case $d=2,$ for each $i = 1,\dots,d-1$ it
is the case that $V(\frak p)$ is the orthogonal direct sum of
$\ker X_i$ and the image of $Y_i.$ Therefore, $V(\frak p)$ is the
orthogonal direct sum of $V(\frak p;0)$ and images of the various
$Y_i.$ If $\frak p = (n_1,\dots,n_d)$ and
$\frak q = (n_1,\dots,n_{i-1}-1,n_i+1,\dots,n_d)$ are both  non--
increasing partitions of $n$ into $d$ parts differing only in
places $i-1$ and $i$ then we write $Y_i\frak q = \frak p.$ The
case $d=2$ shows that $Y_i|V(\frak q;0) \to V(\frak p;0)$ is then
an isomorphism. Therefore no simple modules can appear in the
orthogonal complement of $V(\frak p;0)$ inside $V(\frak p)$
except those which are images of simple modules inside the
various $V(\frak q)$ for which there is a $Y_i$ with $Y_i\frak q
= \frak p.$  Every simple module must therefore come from some
$V(\frak p;0).$ This shows moreover that $V^n$ is a direct sum of
simple submodules; it is ``semisimple'' or ``completely
reducible''. To show that $\Cal H_n$ is a direct sum of total
matric algebras it is therefore sufficient to show that for $d
\ge n$ its representation on $V^n$ is faithful, which is true in
the classical case and therefore also true for generic $q.$ Now
the ranks of the simple $V(\frak p;0)$ modules do not depend on
$q,$ so for $d=n$ we already know that $\Cal H_n$ has a
homomorphic image with total rank equal to that of $\Cal H_n$
itself, namely, $n!.$ Over a field we would be done. As our
coefficient ring is a domain, it follows that the kernel contains
only torsion elements of which there are none since $\Cal H_n$ is
free. $\quad\square$ \enddemo

We have thus come to one of our main conclusions: The Hecke
algebra $\Cal H_n$ ``splits'' over $\Bbb Z[q,q^{-1},1/n_{q^2}!],$
i.e., over that ring of coefficients, it becomes a direct sum of
total matric algebras. This much could already have been obtained
by repairing the omissions in \cite{DJ} or \cite{W}, but the
association to each non-increasing partition of $n$ of a
canonical simple module with an inner product and canonical
orthogonal basis indexed by Young diagrams is not immediately
deducible from either. Since our ``canonical'' orthogonal basis
is indeed a basis, we have the following number--theoretic
assertion.
\proclaim{Proposition 4.3} (i) The norms of all canonical basis
elements of the canonical $\Cal H_n$ modules are units in $\Bbb
Z[q,q^{-1}, 1/n_{q^2}!];$ (ii) (The case $q = 1$) The norms of
all canonical basis elements of the canonical $S_n$ modules are
units in $Z[1/n!]$. $\quad\square$ \endproclaim

\subhead 5. The Donald-Flanigan conjecture \endsubhead

Consider $\Cal H_n = \Cal H_n(q)$ for the moment as defined over
$\Bbb Z_q =\Bbb Z[q,q^{-1}]$ and replace $q$ by $1+t.$  Then
$\Bbb Z_q \subset \Bbb Z[[t]],$ so we can consider $\Cal
H_n(1+t)$ as an algebra over the latter, and as such it is
obviously a deformation of $\Bbb ZS_n.$
\proclaim {Corollary 1.2} Setting $q = 1+t,$ the Hecke algebra
$\Cal H_n(1+t)$ together with the multiplicatively closed subset
$S$ of $\Bbb Z[[t]]$ generated by $1/n_{q^2}! = 1/n_{(1+t)^2}! $
is a global solution to the Donald--Flanigan  problem for the
symmetric group $S_n.$
\endproclaim
\demo{Proof} Not only is  $S^{-1}\Cal H_n$ is separable, but we
have seen that it is already a direct sum of total matric
algebras. Now observe that $n_{q^2}!$ with $q$ replaced by $1+t$
becomes an element of $\Bbb Z[[t]]$ which does not vanish after
reduction modulo any prime,  for no factor
$i_{q^2} = (1-q^{2i})/(1-q^2)$ can vanish.  (If $p|i$ then the
constant term of $i_{q^2}!$ vanishes, but never the entire
expression, so it remains invertible as an element of $\Bbb
F_p((t))$ for every rational prime $p.$)  It follows that $S,$
considered as a subset of  $\Bbb Z[[t]],$ contains no rational
prime.
$\quad\square$\enddemo

\subhead 6. Idempotents \endsubhead

Having completed our principal task, we briefly discuss primitive
idempotents, mainly to show that even in the classical case the
``canonical'' idempotents generated by our procedure differ from
the standard Frobenius--Young idempotents constructed using Young
symmetrizers.  Since the quantized case will take us too far
afield, we restrict attention to $\Bbb Z[1/n!]S_n.$  Suppose that
the rank of $V$ is now also precisely $n.$  Then we now know that
the foregoing group ring is just
$\bigoplus \operatorname{End} V(\frak p;0),$ where $\frak p$ runs
over the non-increasing partitions of $n$ and the endomorphisms
are with respect to the coefficient ring
$\Bbb Z[1/n!].$  Each $V(\frak p;0)$ is free with a canonical
orthogonal basis indexed by Young diagrams, and the norm of each
basis element is a unit.  If $v \in V(\frak p;0)$ is such a basis
element, then the projection of an arbitrary $u\in V(\frak p;0)$
is $((u,v)/(v,v))v,$ which is again in $V(\frak p;0)$ since the
denominator is a unit. Therefore, associated to each Young
diagram there is a primitive idempotent of $\Bbb Z[1/n!]S_n,$
namely, the projection on the corresponding basis vector.   We
call this the {\it canonical} idempotent associated to the Young
diagram.

Consider the simplest non-trivial case, $n=3,$ where we denote
the basis vectors of $V$ by $x_1, x_2, x_3.$  We have three
simple modules, $V((3);0), V((1^3);0),$ and $V((2,1);0).$  The
first two each have rank one and are spanned, respectively, by
$x_1^3$ (not  the ``symmetrized" element
$\sum_{\sigma \in S_n}x_{\sigma1}x_{\sigma2}x_{\sigma3},$ which
lies in $V((1^3))$ and by the skew element
$\sum_{\sigma \in S_n}
(-1)^{\sigma}x_{\sigma1}x_{\sigma2}x_{\sigma3}.$  The module
$V((2,1);0)$ has rank two.  Writing $x$ for $x_1$ and $y$ for
$x_2,$ it has, as we have seen, the canonical basis
$$ v_1 = yx^2-\frac 12 x(xy+yx), \qquad  v_2 = x(xy-yx)$$
corresponding to the Young diagrams which we may write, in
obvious notation, as $[[1,2],[3]]$ and $[[1,3],[2]].$ The
classical Frobenius--Young idempotents corresponding to these
are, respectively,
$\frac 13(1-(13))(1+(12)) = \frac 13(1+(12)-(13)-(123)),$ and
$\frac 13(1-(12))(1+(13)) = \frac 13(1+(13)-(12)-(132)).$
Their sum is the central idempotent $\frac 13(2-(123)-(132)),$
which is the projection operator of $\bigoplus V(\frak p;0)$ on
$V((2,1);0).$  This depends only on the partition $(2,1)$ of $3,$
so we may denote it $e_{(2,1)}.$  Then the projection on $v_1$ is
$\frac 12(1+(23))e_{(2,1)}$ and that on $v_2$ is
$\frac 12(1-(23))e_{(2,1)},$ as one can immediately verify.
These are the ``canonical'' primitive idempotents associated to
$[[12],3]$ and $[[13],2]$ and are clearly not the same as the
Frobenius--Young idempotents (although they necessarily have the
same sum as the corresponding Frobenius--Young idempotents).

Notice that once we know the central idempotents of $\Bbb
Z[1/n!]S_n$ then the canonical idempotents can be computed
inductively from the case of $n-1$ (with $S_{n-1}$ acting on the
last $n-1$ variables) because of the way that $V(\frak p;0)$ is
constructed from the $V(\frak p';0)$ with $\frak p' \to \frak p.$
 This is illustrated here.  The basis element $v_1$ has come from
the unique canonical basis element $xx$ of $V((2);0)$ and $v_2$
from the unique basis element $xy-yx$ of $V((1,1);0).$  The
projections on these are $\frac12(1+(12))$ and $\frac12(1-(12)),$
respectively.  These give rise to the projections on $v_1$ and
$v_2$ where $(12)$ must now be replaced by $(23)$ since $S_2$ is
now operating on the last two variables.  The central idempotents
can be computed using the independent Casimir operators of
degrees $2,\dots,n$ of $Usl_n.$  This is not difficult since we
only need to know how they operate on the $V(\frak p;0),$ which
is given readily by a beautiful theorem of Harish--Chandra (cf.
\cite{K, p.308}).  Details, particularly in the quantized case,
are the subject for another paper.

\subhead 7. Afterword: Status of the Donald--Flanigan conjecture
and applications to group theory \endsubhead

The general Donald--Flanigan conjecture, even in its original
form, has proven remarkably resilient. It has been
characterized as a ``modular form of Maschke's theorem.'' The
statement is straightforward but seems to go deep into the
structure of finite groups. It encodes certain assertions about
the cohomology of groups which imply that for any finite group
$G$ and prime $p$ dividing $\#G$ one has $H^1(G, \Bbb F_p) \ne
0$. This in turn implies a previously unknown ``dual'' to
Cauchy's theorem: If $p|\#G$ then there is an element $g \in G$
whose centralizer $C_G(g)$ contains a normal subgroup of index
$p$, cf. \cite{GG}.  This has been verified by Fleischmann,
Janiszczak and Lempken \cite{FJL} who prove the stronger ``weak
non--Schur property'': There  exists in  $G$  an element  $g$
whose order is divisible by  $p$  and  whose centralizer $C =
C_G(g)$ has the property that its  commutator subgroup $C'$ does
not contain the $p$--part of $g$.   This proposition can be
reduced to the case of  simple  $G$ and they prove it using the
most difficult result so far known in group theory -- the
classification of the finite simple groups -- by showing that it
is  indeed true for all of them.  There is some delicacy to the
choice of hypothesis; the  ``strong non-Schur property''
which  asserts that the preceding  $g$  can be taken to be a
$p$--element  fails for three of the exceptional groups. It is an
interesting  and not too difficult exercise to prove the weak
non--Schur  property directly for the alternating and symmetric
groups.  The existence of a global solution may have even more
cohomological and group--theoretic consequences.

Various cases of the Donald--Flanigan conjecture are known, but
so far, no other global equimodular solutions except the
one in \cite{ES}. Donald and Flanigan in their original 1974
paper \cite{DF} settled only the case of abelian groups. There
was no further  progress until 1988, when the second author
\cite{Sps} proved the conjecture for the group algebra over $\Bbb
F_p$ of groups with cyclic $p$--Sylow subgroups (a condition
equivalent, by Higman's theorem, to the group algebra having
finite representation type). The present authors \cite{GSps} also
settled the  case of  groups with abelian normal $p$--Sylow. The
Donald--Flanigan conjecture has been verified for all groups of
order less than 32 except for the ``extra special" group of order
27. In \cite{Sps2}, the second author showed that for $p$--
solvable  groups with  cyclic $p$--Sylow, the integral group ring
and the    $p$--modular  semisimple  deformation can be achieved
by a single deformation with a ``discrete'' and a  ``continuous
axis''. Such a deformation has been called {\it liftable}.
Michler \cite{Mi} proved a local version of this  result for
blocks of  cyclic defect group.   In \cite{ES}, the conjecture is
confirmed for all blocks with   dihedral  defect group, in
\cite{Sps3}, the deformation of the dihedral  2-groups are shown
to be liftable, and in \cite{ESps} this result is extended to
other dihedral groups.  One might expect that the Hecke algebra
of a finite Coxeter group is always a global solution to the
Donald--Flanigan problem but this fails, e.g., for the dihedral
group $D_n$ of order $2n.$ In that case, however, there is a
suitable deformation of the Hecke algebra which serves,
suggesting that for finite Coxeter groups there is always a
global solution which is at least a deformation of its Hecke
algebra.

While there is thus considerable evidence for the truth of the
Donald--Flanigan conjecture, it seems to lie deep and remains one
of the most tantalizing open problems in finite group theory.

\enddocument